\documentstyle[12pt]{article}

\newcommand{\qed}{$\;\;\;\Box$}
\newenvironment{proof}{\par\smallbreak{\sl Proof.~}}
{\unskip\nobreak\hfill \qed \par\medbreak}

\newtheorem{thm}{Theorem}

\newtheorem{prop}[thm]{Proposition}
\newtheorem{defn}[thm]{Definition}
\newtheorem{cor}[thm]{Corollary}
\newtheorem{rem}[thm]{Remark}
\newtheorem{ex}[thm]{Example}

\title{
Initial-Boundary Problems for Semilinear Hyperbolic Systems
with Singular Coefficients
}

\newcounter{thesame}
\author{
I.~Kmit\\
{\small
Institute for Applied Problems of Mechanics and Mathematics,}\\
{\small Ukrainian Academy of Sciences}\\
{\small Naukova St.\ 3b}\\
{\small 79060 Lviv,
Ukraine}
\\
{\small   E-mail:
{\tt kmit@ov.litech.net}}
}

\date{}

\begin{document}

\maketitle

\begin{abstract}
In the paper we use the framework of Colombeau algebras of generalized
functions to study existence and uniqueness  of global
generalized solutions to  mixed non-local problems for a semilinear
hyperbolic system. Coefficients of the system
 as well as initial and boundary data are allowed to be strongly
singular, as the Dirac delta function and derivatives thereof.
To obtain the existence-uniqueness result we prove a criterion of
invertibility in the full version of the Colombeau algebras.

\end{abstract}

\section{Introduction}\label{sec:intr}
In the domain $\Pi =\{(x,t)|-L<x<L$, $t>0\}$
we consider the following initial-boundary value problem for a
generalized function $U$:
\begin{eqnarray}
(\partial_t  + \Lambda(x,t)\partial_x) U &=& f(x,t,U), \qquad (x,t)\in \Pi
  \label{eq:1}\\
U|_{t=0} &=& A(x), \qquad x\in (-L,L) \label{eq:2} \\
B(t)U|_{x=-L}+C(t)U|_{x=L} +\int\limits_{-L}^{L}D(x,t)U\,dx
  &=& H(t), \qquad  t\in(0,\infty) \label{eq:3} \, .
\end{eqnarray}
where,  $U$, $f$, $A$, and $H$ are real $n$-vectors, $\Lambda$, $B$,
$C$,  and $D$ are real
$(n\times n)$-matrices, and  $\Lambda=\mbox{\rm diag}(\Lambda_1,\dots,\Lambda_n)$
 is a diagonal
matrix.

Special cases of~(\ref{eq:1})--(\ref{eq:3}) are mathematical
formulations of problems arising in
population dynamics \cite{bo-ar,Huy,son},
laser dynamics \cite{r1,rw,sbw,tlo}, and
chemical kinetics~\cite{z}.

Our goal is to find global solutions to problem~(\ref{eq:1})--(\ref{eq:3})
when the data $\Lambda$, $A$, $B$, $C$, $D$, and $H$ are allowed to be
strongly singular (at least of the Dirac delta type). This entails
multiplication of distributions in~(\ref{eq:1}) and~(\ref{eq:3}). Indeed,
since initial singularities expand from $\partial\Pi$ into $\Pi$ along
characteristic curves of~(\ref{eq:1}), one can expect that solutions
within $\Pi$ are at least as singular as they are on $\partial\Pi$.
Furthermore, since the characteristics of~(\ref{eq:1})  are singular
themselves, we also meet the problem of composition of two
singular functions (for instance, the composition of singular
initial data with singular characteristic curves).
It is known~\cite{Courant} that
even if $F$ is a regular function, but $S$ is a singular one, then
$F(S(x))$ is not well-defined in ${\cal D}'$.
Finally, it should be noted that such three ingredients as
singularities, nonlinear operations, and differentiation, cannot be
presented unrestrictedly within ${\cal D}'$. All this makes impossible
to use the framework of the distribution theory for our purpose.
Nevertheless, such a differential-algebraic structure
as an algebra of generalized functions is able to deal with the above
problems in a quite reasonable way.
We here  use the  Colombeau version ${\cal G}$  of an algebra, which is defined on
any domain in ${\bf R}^n$ as well as on its   closure,
is a sheaf, and admits  restrictions to the coordinate planes.

We hence assume that entries of $A$ are generalized functions
in the Colombeau algebra ${\cal G}[-L,L]$, entries of $B$, $C$, and $H$
are from ${\cal G}({\bf R}_{+})$, and entries of $\Lambda$ and $D$ are
from ${\cal G}(\overline{\Pi})$.

Another advantage of using Colombeau algebra of generalized functions
lies in the fact that in a variety of important cases
the division by generalized functions, in particular the
division by discontinuous functions and measures, is defined in ${\cal G}$.
The latter, of course, is impossible in ${\cal D}'$.
We completely describe the cases when the division is possible by
obtaining a criterion of
invertibility in ${\cal G}(\Omega)$.

The plan of our exposition is as follows. Section 2 presents some
preliminaries.
In Section~3 we extend  the criterion  of invertibility from the
simplified version of Colombeau algebra ${\cal G}_s(\Omega)$ (see~\cite{GKOS})
to its full version ${\cal G}(\Omega)$. The main result of the paper
is given in  Section 4, where we prove  the global
existence-uniqueness theorem within  ${\cal G}(\Omega)$.

A novelty of the paper is that it treats singular coefficients
in~(\ref{eq:1}) in the context of mixed problems for a quite wide
range of boundary conditions which can be classical as well as
nonclassical (nonseparable and integral).

Existence-uniqueness results within Colombeau algebras for two-dimensional
hyperbolic problems with discontinuous coefficients were studied
in~\cite{hoop,LOb,Obe,Obe2}. Note that
the discontinuity implies global
boundedness estimates on the coefficients within Colombeau algebra ${\cal G}$, thereby avoiding
the negative effect of infinite propagation speed.  At the present paper
we do not impose the assumption of global boundedness on coefficients
of~(\ref{eq:1}), thereby allowing them to be strongly
singular. In~\cite{NP} the authors use the Colombeau algebra of
tempered generalized functions ${\cal G}_{\tau}$ to succeed with strongly
singular coefficients in Cauchy problems for hyperbolic systems.

\section{Preliminaries}

In this section we summarize the relevant
material on Colombeau algebras of generalized functions.

Let $\Omega\subset{\bf R}^n$ be a domain in ${\bf R}^n$. We denote by ${\cal G}(\Omega)$
and ${\cal G}(\overline\Omega)$ the full version of Colombeau algebra of
generalized functions over $\Omega$ and $\overline\Omega$, respectively.
To define ${\cal G}(\Omega)$ and ${\cal G}(\overline\Omega)$, we
introduce the mollifier spaces in order to parametrize
the regularizing sequences of generalized functions.
For $q\in {\bf N}_0$ denote
$$
\begin{array}{cc}
\displaystyle
{\cal A}_q({\bf R})=\Bigl\{\varphi\in{\cal D}({\bf R})\Bigm|\int\varphi(x)\,dx=1,\int x^k\varphi(x)\,dx=0
\mbox{\ for \ } 1\le k\le q\Bigr\},\\[7mm]
{\cal A}_q({\bf R}^n)=\Bigl\{\varphi(x_1,\dots,x_n)=\prod\limits_{i=1}^n\varphi_0(x_i)\Bigm|
\varphi_0\in{\cal A}_q({\bf R})\Bigr\}.
\end{array}
$$
For $\varphi\in{\cal A}_0({\bf R}^n)$ define
$$
\varphi_{\varepsilon}(x)=\frac{1}{\varepsilon^n}\varphi\biggl(
\frac{x}{\varepsilon}\biggr).
$$
We now introduce the algebra of moderate elements
${\cal E}_{\mbox{\rm M}}(\overline{\Omega})$ in the
following way.  Define
$$
{\cal E}(\overline{\Omega})=\{u: {\cal A}_0\times\overline{\Omega}\to {\bf R} \Bigm|
u(\varphi,.)\in C^{\infty}(\overline{\Omega})\ \  \forall \varphi\in{\cal A}_0({\bf R})\}.
$$
Now  ${\cal E}_{\mbox{\rm M}}(\overline{\Omega})$, is defined to be a
subalgebra of ${\cal E}(\overline{\Omega})$  consisting of
elements $u\in{\cal E}(\overline{\Omega})$ with the
following property:
$$
\begin{array}{cc}
\forall K\subset \overline{\Omega} \,\,\mbox{\rm compact \ }, \forall \alpha\in{\bf N}_0^n,
\exists N\in{\bf N}\,\, \mbox{\rm such\  that \ } \forall \varphi\in{\cal A}_N({\bf R}^n),
 \\[5mm]
\displaystyle
\exists C>0, \exists \eta>0\,\,
\mbox{\rm with \ }
\sup\limits_{x\in K}|\partial^{\alpha}u(\varphi_{\varepsilon},x)|\le C
\varepsilon^{-N}, \ \ 0<\varepsilon<\eta.
\end{array}
$$
The ideal ${\cal N}(\overline{\Omega})$ consists of all $u\in{\cal E}_{\mbox{
\rm M}}(\overline{\Omega})$
such that
$$
\begin{array}{cc}
\forall K\subset \overline{\Omega} \,\,\mbox{\rm compact \ }, \forall \alpha\in{\bf N}_0^n,
\exists N\in{\bf N} \,\,\mbox{\rm such\  that \ } \forall q\ge N,
\forall \varphi\in{\cal A}_q({\bf R}^n),
 \\[5mm]
\displaystyle
\exists C>0, \exists \eta>0
\,\,\mbox{\rm with \ }
\sup\limits_{x\in K}|\partial^{\alpha}u(\varphi_{\varepsilon},x)|\le C\varepsilon^{q-N}, \ \ 0<\varepsilon<\eta.
\end{array}
$$
Finally,
$$
{\cal G}(\overline{\Omega})={\cal E}_{\mbox{\rm M}}(\overline{\Omega})/{\cal N}(\overline{\Omega}).
$$
This is an associative, commutative differential algebra.
The algebra ${\cal G}({\Omega})$ on open set is constructed in the same manner
(with $\Omega$ in place of $\overline{\Omega}$ in the definition above). Note that
${\cal G}({\Omega})$ admits a
canonical embedding of ${\cal D}'(\Omega)$. We will use the notation
$U=[(u(\varphi,x))_{\varphi\in{\cal A}_0({\bf R}^n)}]$ for elements $U$ of
Colombeau algebra ${\cal G}(\Omega)$ with $u(\varphi,x)$ to be a representative of
$U$.

To reduce information from generalized functions to the level of
distributions, we use the notion of an associated distribution. We
say that $U\in{\cal G}(\Omega)$ admits
$f\in{\cal D}'(\Omega)$ as associated distribution
(or $U$ is associated to $f$), denoted by $U\approx f$, if for all
$\psi\in{\cal D}(\Omega)$ there exists $N\in{\bf N}$ such that
$$
\lim\limits_{\epsilon\to 0}\int u(\varphi_{\varepsilon},x)\psi(x)\,dx=\langle
f,\psi\rangle
$$
for all $\varphi\in{\cal A}_N({\bf R}^n)$.

\section{Criterion of invertibility in the full version of
\mbox{} Colombeau algebra of generalized functions}

In spite of the fact that ${\cal G}(\Omega)$ is not a field, the division by
singular distributions (in particular, by discontinuous functions and
measures) is sometimes possible. It is given by the
following criterion of (multiplicative) invertibility.

Let $\Omega\subset{\bf R}^n$ be an arbitrary subdomain in ${\bf R}^n$.

\begin{thm}\label{thm:invert}
Let $U\in{\cal G}(\Omega)$ ($U\in{\cal G}(\overline\Omega)$). Then the following two
conditions are equivalent:\\
(i) $U$ is invertible in ${\cal G}(\Omega)$
(in ${\cal G}(\overline\Omega)$), i.e., there exists
$V\in{\cal G}(\Omega)$ ($V\in{\cal G}(\overline\Omega)$) such that $UV=1$ in ${\cal G}(\Omega)$
(in ${\cal G}(\overline\Omega)$).\\
(ii) For each representative $(u(\varphi,x))_{\varphi\in{\cal A}_0({\bf R}^n)}$ of $U$
and each compact set $K\subset\Omega$ ($K\subset\overline\Omega$)
there exists $p\in{\bf N}$ such that for all $\varphi\in{\cal A}_p({\bf R}^n)$ there is
$\eta>0$ with $\inf\limits_{K}|u(\varphi_{\varepsilon},x)|
\ge\varepsilon^p$ for all $0<\varepsilon<\eta$.
\end{thm}

Note that the criterion of invertibility for the simplified version of Colombeau
algebra ${\cal G}_s(\Omega)$, where $\Omega$ is open, was proved in~\cite{GKOS}.

\begin{proof}
We use the argument similar to that presented in~\cite{GKOS}. We prove
the desired assertion for an arbitrary fixed open set $\Omega$ (the proof
for the closed set
$\overline\Omega$ is similar).

$(i)\Rightarrow(ii)$. Set $U=[(u(\varphi,x))_{\varphi\in{\cal A}_0({\bf R}^n)}]$ and
$V=[(v(\varphi,x))_{\varphi\in{\cal A}_0({\bf R}^n)}]$. By assumption, there exists
$N=[(n(\varphi,x))_{\varphi\in{\cal A}_0({\bf R}^n)}]\in{\cal N}(\Omega)$ such that
$u(\varphi_{\varepsilon},x)v(\varphi_{\varepsilon},x)=1+n(\varphi_{\varepsilon},x)$
for all $\varphi\in{\cal A}_0({\bf R}^n)$.

Fix an arbitrary compact set $K\subset\Omega$. We first prove that there
exists $p\in{\bf N}$ such that for all
$\varphi\in{\cal A}_p({\bf R}^n)$ there is $\eta>0$ with $v(\varphi_{\varepsilon},x)\ne 0$
for all $x\in K$ and $0<\varepsilon<\eta$. Assume, to the contrary, that the latter
is not true. This means that for each $p\in{\bf N}$ there exist
$\varphi\in{\cal A}_p({\bf R}^n)$, a sequence $\varepsilon_n\searrow 0$, and a
sequence $x_n\in K$ such that $v(\varphi_{\varepsilon_n},x_n)=0$ for all $n\ge 1$.
Hence $0=u(\varphi_{\varepsilon_n},x_n)v(\varphi_{\varepsilon_n},x_n)=
1+n(\varphi_{\varepsilon_n},x_n)$ and finally $n(\varphi_{\varepsilon_n},x_n)=-1$
for all $n\ge 1$, a contradiction to the fact that $N\in{\cal N}(\Omega)$.
Since $v\in{\cal E}_{\mbox{\rm M}}(\Omega)$, there exists $q\in{\bf N}$ such that for all
$\varphi\in{\cal A}_q({\bf R}^n)$ there are $C>0$ and $\mu>0$ with
$\sup\limits_{K}|v(\varphi_{\varepsilon},x)|\le C/\varepsilon^q$ for all
$0<\varepsilon<\mu$. Set $\tilde q=max\{p,q\}$.
Due to the fact that ${\cal A}_{q+1}({\bf R}^n)\subset{\cal A}_q({\bf R}^n)$ for all $q\in{\bf N}_0$,
we conclude that for each
$\varphi\in{\cal A}_{\tilde q}({\bf R}^n)$
the estimate
$$
\inf\limits_{K}|u(\varphi_{\varepsilon},x)|\ge
\frac{\varepsilon^q}{C}\Bigl(1-\sup\limits_K|n(\varphi_{\varepsilon},x)|
\Bigr)\ge
\varepsilon^{q+1}
$$
is true for all sufficiently small $\varepsilon$.
Since $K$ is an arbitrary compact subset of $\Omega$, the desired assertion
follows.

$(ii)\Rightarrow(i)$. Consider a covering $(K_i)_{i\in{\bf N}}$ of $\Omega$
by compact sets
$K_i$ such that $K_1\subset K_2\subset\dots\subset\Omega$. It is known that,
if $W\in{\cal G}(K_{i+1})$, then $W|_{K_j}\in{\cal G}(K_{j})$ for all $j\le i$. This
fact is true due to the sheaf properties of ${\cal G}(\Omega)$.

Set $v(\varphi_{\varepsilon},x)=1/u(\varphi_{\varepsilon},x)$
and $v_i(\varphi_{\varepsilon},x)=v(\varphi_{\varepsilon},x)|_{K_i}$.
Fix an arbitrary $i\in{\bf N}$. By assumption, there exists $p\in{\bf N}$ such that for all
$\varphi\in{\cal A}_p({\bf R}^n)$ there is a constant $\eta(\varphi)>0$
such that the expression
$1/u(\varphi_{\varepsilon},x)$
exists for all $0<\varepsilon<\eta(\varphi)$ and for all $x\in K_i$. For each
$\varphi\in{\cal A}_p({\bf R}^n)$ let us set
$v_i(\varphi_{\varepsilon},x)
\equiv 0$, where $0<\varepsilon<\eta(\varphi)$ and $x\in K_i$.
Consider the map $\varphi\to
v_i(\varphi,x)\,:\,{\cal A}_0({\bf R}^n)\to C^{\infty}(K_i)$. Let us show that
this map is moderate. Indeed, for each $\varphi\in{\cal A}_p({\bf R}^n)$ we have
$$
\sup\limits_{K_i}|v_i(\varphi_{\varepsilon},x)|=
\frac{1}{\inf\limits_{K_i}|u(\varphi_{\varepsilon},x)|}
\le\frac{1}{\varepsilon^{p}}
$$
for all sufficiently small $\varepsilon>0$. The moderate estimate
for $\partial^{\alpha}
v_i(\varphi,x)$, where $|\alpha|=1$, follows from the simple estimate
$$
\biggl|\partial^{\alpha}\biggl(\frac{1}{u(\varphi_{\varepsilon},x)}
\biggr)\biggr|=
\biggl|\frac{\partial^{\alpha}(u(\varphi_{\varepsilon},x))}{u^2(\varphi_{\varepsilon},x)}
\biggr|\le\frac{\partial^{\alpha}(u(\varphi_{\varepsilon},x))}{\varepsilon^{2p}}
$$
for sufficiently small $\varepsilon$
and from
the moderateness of
$\partial^{\alpha}u(\varphi,x)$. Proceeding similarly with the
higher-order derivatives of $v_i$, we conclude that
$[(v_i(\varphi,x))_{\varphi\in{\cal A}_0({\bf R}^n)}]$, denoted by
$V_i$, belongs to
${\cal E}_{\mbox{\rm M}}(K_i)$. Furthermore, it is the inverse to $U$
in ${\cal G}(K_i)$.
{}From the definition of Colombeau generalized functions and the construction
of $V_i$ it follows that $V_i|_{K_j}\in{\cal G}(K_j)$ for all $j\le i$. We
therefore obtained a coherent family $\{V_i,i\in{\bf N}\}$. By the sheaf properties
of ${\cal G}(\Omega)$, there exists a unique element $V\in{\cal G}(\Omega)$ such that
$V|_{K_i}\in{\cal G}(K_i)$ for all $i\ge 1$. By construction, $V$ is an inverse to
$U$ in ${\cal G}(\Omega)$.
\end{proof}

We now take into account the definition of Colombeau generalized numbers
and the fact that an element $U\in{\cal G}(\Omega)$ is a constant iff there is
$r\in\overline{\cal C}$ such that $U-r=0$ in ${\cal G}(\Omega)$. The
following corollary provides a criterion of invertibility of Colombeau
generalized numbers within the full version of Colombeau algebras.
For the same result within the simplified version of Colombeau algebras
see~\cite{GKOS}.

\begin{cor}\label{cor:invert}
Let $r\in\overline{\cal C}$. Then the following two conditions are
equivalent:\\
(i) $r$ is invertible in $\overline{\cal C}$, i.e., there exists
$s\in\overline{\cal C}$ with $rs=1$ in $\overline{\cal C}$.
\\
(ii) For each representative
$(r(\varphi))_{\varphi\in{\cal A}_0({\bf R})}$
of $r$ there exists $p\in{\bf N}$ such that for all $\varphi\in{\cal A}_p({\bf R})$ there is
$\eta>0$ with $|r(\varphi_{\varepsilon})|\ge\varepsilon^p$ for all $0<\varepsilon<\eta$.
\end{cor}

\begin{ex}\label{ex:delta}\rm
Let
$$
U=\biggl[\biggl(l(\varphi)+\frac{1}
{l^{m+1}(\varphi)}\Phi^{(m)}\biggl(\frac{x}
{l(\varphi)}\biggr)\biggr)_{\varphi\in{\cal A}_0({\bf R})}\biggr]\in{\cal G}(\Omega),
$$
where $\Omega\subset{\bf R}$,
$l(\varphi)=\sup\{|y|,\varphi(y)\ne 0\}$, $\Phi(x)\in{\cal D}(\Omega)$
is a fixed element of ${\cal D}(\Omega)$ such that $\int\Phi(x)\,dx=1$ and
$\Phi(x)\ge 0$.
One can easily see that
$U\approx\delta^{(m)}$. Indeed, for an arbitrary $\psi(x)\in{\cal D}(\Omega)$ we have
$$
\lim\limits_{\varepsilon\to 0}\int \biggl(\varepsilon l(\varphi)+\frac{1}
{\varepsilon l^{m+1}(\varphi)}\Phi^{(m)}\biggl(\frac{x}
{\varepsilon l(\varphi)}\biggr)\biggr)\psi(x)\,dx=<\delta^{(m)},\psi>.
$$
Since $l(\varphi_{\varepsilon})=\varepsilon l(\varphi)$, we have the
following estimate: for each compact set $K\subset\Omega$ and for each
$\varphi\in{\cal A}_2({\bf R})$ there exists $\eta>0$ with
$$
\inf\limits_{K}\biggl|l(\varphi)+\frac{1}
{l^{m+1}(\varphi)}\Phi^{(m)}\biggl(\frac{x}
{l(\varphi)}\biggr)\biggr|=\inf\limits_{K}\biggl|\varepsilon l(\varphi)+
\frac{1}{\varepsilon l^{m+1}(\varphi)}\Phi^{(m)}
\biggl(\frac{x}
{\varepsilon l(\varphi)}\biggr)\biggr|\ge\varepsilon^2,\quad 0<\varepsilon<\eta.
$$
This estimate  is uniform with respect to all compact sets $K\subset\Omega$
and $\varphi\in{\cal A}_2({\bf R})$. By Theorem~\ref{thm:invert},
$U$ is invertible in ${\cal G}(\Omega)$.
\end{ex}
This example shows that within ${\cal G}(\Omega)$ the division by the derivatives of
the delta-function is possible.

\begin{prop}\label{prop:unique}
Let $U\in{\cal G}(\Omega)$ ($U\in{\cal  G}(\overline\Omega)$)
and $U$ is invertible in ${\cal G}(\Omega)$ (in ${\cal G}(\overline\Omega)$). Then the
multiplicative inverse of $U$ is unique.
\end{prop}
\begin{proof}
We prove
the desired assertion for an open set $\Omega$ (the proof for the closed set
$\overline\Omega$ is similar).

Assume, to the contrary, that $U$ possesses two multiplicative inverses
$V_1,V_2\in{\cal G}(\Omega)$. This implies the equality
$$
U(V_1-V_2)=0\,\,\mbox{in}\,\,{\cal G}(\Omega).
$$
We conclude from
Theorem~\ref{thm:invert}, specifically from the local invertibility
estimate, that $U\not\in{\cal N}(\Omega)$,
hence that $V_1-V_2\in{\cal N}(\Omega)$, and finally that $V_1=V_2$ in
${\cal G}(\Omega)$, a contradiction to our assumption.
\end{proof}

\section{Existence and uniqueness of Colombeau generalized solutions}

In this section we develop the results of~\cite{HKm} to the case
of singular coefficients in~(\ref{eq:1}).
Simultaniously, we consider less restrictive conditions on the initial data
in~(\ref{eq:2}) and~(\ref{eq:3}).
To prove a general global existence and uniqueness result in Colombeau
algebra of generalized functions, we need the following definition of
generalized functions of a less restrictive growth if comparing with
$1/\varepsilon$-growth (see the definition of ${\cal E}_{\mbox{\rm M}}$).

\begin{defn}\label{defn:ga}\rm
Let $\Omega\subset{\bf R}^n$ be a domain in
${\bf R}^n$.
Suppose we have a  function $\gamma : (0,1)\mapsto (0,\infty)$.
An element $U\in{\cal G}(\Omega)$ ($U\in
{\cal G}(\overline{\Omega})$) is called
{\it locally of $\gamma$-growth}, if it has a representative
$u\in{\cal E}_{\mbox{\rm M}}(\Omega)$
($u\in{\cal E}_{\mbox{\rm M}}(\overline\Omega)$) with the following property:

For every compact subset $K\subset\Omega$
there is $N\in{\bf N}$ such that for every $\varphi\in {\cal A}_{{\cal N}}({\bf R}^n)$
there exist $C>0$ and $\eta>0$ with
$\sup\limits_{x\in K}|u(\varphi_{\varepsilon},x)|\le C
\gamma^N(\varepsilon)$ for
$0<\varepsilon<\eta$.

\end{defn}
Note that this definition generalizes Definition~2 from~\cite{HKm}.

\begin{defn}\label{defn:gainvert}\rm
Let $\Omega\subset{\bf R}^n$ be a domain in
${\bf R}^n$.
Suppose we have a function $\gamma : (0,1)\mapsto (0,\infty)$.
An element $U\in{\cal G}(\Omega)$ ($U\in
{\cal G}(\overline{\Omega})$) is called
{\it locally $\gamma$-invertible}, if it has a representative
 $(u(\varphi,x))_{\varphi\in{\cal A}_0({\bf R}^n)}$
with the following property:

for each compact set $K\subset\Omega$ ($K\subset\overline\Omega$)
there exists $p\in{\bf N}$ such that for all $\varphi\in{\cal A}_p({\bf R}^n)$ there is
$\eta>0$ with $\inf\limits_{K}
|u(\varphi_{\varepsilon},x)|\ge\gamma^{-p}(\varepsilon)$ for all
$0<\varepsilon<\eta$.

\end{defn}


We now make several assumptions on the initial data
of problem~(\ref{eq:1})--(\ref{eq:3}). Let $\gamma(\varepsilon)$ and
$\gamma_1(\varepsilon)$ be positive functions from
$(0,1)$ to $(0,\infty)$ having the properties
\begin{equation}\label{eq:ga}
{\gamma(\varepsilon)}^{\gamma^N(\varepsilon)}=
O\biggl(\frac{1}{\varepsilon}\biggr),
\quad {\gamma_1(\varepsilon)}^{\gamma_1^N
(\varepsilon)}
=O\biggl(\frac{1}{\varepsilon}\biggr),
\quad {\gamma(\varepsilon)}^{\gamma_1^N
(\varepsilon)}
=O\biggl(\frac{1}{\varepsilon}\biggr)\quad\mbox{as}\,\,\varepsilon\to 0
\end{equation}
for each $N\in{\bf N}$. Assume that

\begin{enumerate}
\item The mapping $U\mapsto f(x,t,U)$ and all its derivatives are polynomially
  bounded, uniformly over $(x,t)$ varying in compact subsets of
$\overline{\Pi}$.
\item The mapping $U\mapsto \nabla_U f(x,t,U)$ is globally bounded,
  uniformly over $(x,t)$ varying in compact subsets of $\overline{\Pi}$.
\item $\Lambda_1,\dots,\Lambda_k<0$, $\Lambda_{k+1},\dots,\Lambda_n>0$
(these inequalities are satisfied on the level of representatives),
where $k$ is fixed and $1\le k\le n$.

\item $\Lambda_i$  and $D_{ij}$ for $i\le n$
and $j\le n$ are locally of $\gamma$-growth on $\overline\Pi$.
\item $B_{ij}$ and $C_{ij}$ for $i\le n$
and $j\le n$ are locally of $\gamma$-growth on $[0,\infty)$.
\item $\partial_x\Lambda_i$ for $i\le n$
 are locally of $\gamma_1$-growth on $\overline\Pi$.
\item $\Lambda_i$ for $i\le n$ are locally $\gamma$-invertible
on $\overline\Pi$.
\item   The determinant of the matrix
  $$
  R(t)=\left(
  \begin{array}{cccccc}
    B_{1,k+1}&\ldots & B_{1n} & C_{11}&\ldots & C_{1k}\\
    B_{2,k+1}&\ldots & B_{2n} & C_{21}&\ldots & C_{2k}\\
    \vdots&\ddots&\vdots&\vdots&\ddots&\vdots\\
    B_{n,k+1}&\ldots & B_{nn} & C_{n1}&\ldots & C_{nk}
  \end{array}
  \right)
  $$
is locally $\gamma$-invertible
 on $[0,\infty)$.
\item $\mbox{\rm supp} A_i(x)\subset (-L,L)$; $\mbox{\rm supp} B_{ij}(t)$, $\mbox{\rm supp}
  C_{is}(t)\subset(0,\infty)$ for $1\le i\le n$, $1\le j\le k$, $k+1\le s\le
  n$; $\mbox{\rm supp} D_{im}(x,t)\subset(0,\infty)\times[-L,L]$ for $1\le i,m\le n$.
\end{enumerate}

Let $U\in{\cal G}(\Omega)$ and a smooth function $g(x)$ be slowly increasing at the
infinity. By the definition of ${\cal G}(\Omega)$,
we have $g(U)\in{\cal G}(\Omega)$. Due to this fact and Assumption~1, $f(x,t,U)$
is a well-defined element of ${\cal G}(\Omega)$. Condition~2 is, in fact, sufficient
and is imposed to ensure the global classical solvability of
problem~(\ref{eq:1})--(\ref{eq:3})
with smooth initial data.
We need Assumption~7 to transform the initial problem into an
equivalent integral-operator form. Assumption~8 ensures the compatibility of
(\ref{eq:2}) and (\ref{eq:3}) of any desired order. The hyperbolicity of
system~(\ref{eq:1}) is ensured by Assumption~3.

The point of Assumption~4 is that it allows us to consider
$\Lambda_i$, $B_{ij}(t)$, $C_{ij}(t)$, and $D_{ij}(x,t)$  being
discontinuous functions,
the delta functions, and the derivatives thereof.
An illustration of this fact is given by  Example~\ref{ex:delta} if one takes
$\gamma(l(\varphi))$ in place of $1/l(\varphi)$, where $\gamma$
is specified by
$\gamma(p)=\sqrt{\frac{1}{2}\log{\log{\log{(1/p)}}}}$. If one takes in
addition $\gamma_1(p)={\sqrt{\log{\log{(1/p)}}}}$, the same example shows
that
Assumptions~4 and~5 on $\Lambda$ do not
contradict one another.

We are prepared to state the main result of the paper.
\begin{thm}\label{thm:gen}
Suppose that
$A\in{\cal G}[-L,L]$, $\Lambda,D\in{\cal G}(\overline{\Pi})$, $B$, $C$, $H\in {\cal G}({\bf R}_{+})$,
and $f$ is smooth with respect to all its arguments.
Under Assumptions 1--8 where the functions $\gamma$ and $\gamma_1$ are
specified by~(\ref{eq:ga}),
problem~(\ref{eq:1})--(\ref{eq:3})
has a unique solution $U\in{\cal G}(\overline{\Pi})$.
\end{thm}

\begin{proof}
We first transform problem~(\ref{eq:1})--(\ref{eq:3}) into an
equivalent integral-operator form. Note that all algebraic
operations as well as operation of integration over finite
intervals will be carried out on the level of representatives.
Denote by $\omega_i(\tau;x,t)$ the $i$-th characteristic of~(\ref{eq:1})
passing through a point
$(x,t)\in\overline{\Pi}$, i.e., $\xi=\omega_i(\tau;x,t)$ is the solution to
the Cauchy problem:
$$
\frac{d\xi}{d\tau}=\Lambda_i(\xi(\tau),\tau), \quad \xi(t)=x.
$$
The smallest value of
$\tau\geq 0$ at which the characteristic
$\xi=\omega_i(\tau;x,t)$ intersects $\partial\Pi$
will be denoted by $t_i(x,t)$.

By Assumption 7 and Theorem~\ref{thm:invert},
$\det R(t)$ has an inverse  with entries in ${\cal G}(\bar{\Pi})$.
Using in addition  Proposition~\ref{prop:unique}, we conclude that
there exists a unique element $(\det R)^{-1}\in {\cal G}(\overline{\Pi})$ such that
$\det R\;(\det R)^{-1}=1$. This means that the local part of boundary
conditions~(\ref{eq:3})  is solvable with respect to those components of $U$
whose characteristics move into $\Pi$. Using this fact and integrating
each equation of~(\ref{eq:1}) along the corresponding
characteristic curve,
we obtain the following integral-operator
form of~(\ref{eq:1})--(\ref{eq:3}):
\begin{equation}\label{eq:integral}
\begin{array}{cc}
\displaystyle
U_i(x,t)=(R_iU)(x,t)+
\int\limits_{t_i(x,t)}^t\Bigl[U(\omega_i(\tau;x,t),\tau)\int\limits_0^1\nabla_U
f_i(\omega_i(\tau;x,t),\tau,\sigma U)\,d\sigma
\\[5mm]\displaystyle
+f_i(\omega_i(\tau;x,t),\tau,0)\Bigr]\,d\tau,\ \
1\le i\le n,
\end{array}
\end{equation}
where
$$
(R_iU)(x,t)=
\cases{M_i(t_i(x,t))&if
$t_i(x,t)>0$,
 \cr
A_{i}(\omega_i(0;x,t)) &if
$t_i(x,t)=0$, \cr}
$$
$$
\begin{array}{cc}
M_i(t)=U_i|_{x=-L},\ \ k+1\le i\le n\\
M_i(t)=U_i|_{x=L},\ \ 1\le i\le k,
\end{array}
$$
and
$$
\begin{array}{cc}
\displaystyle
M_i(t)=\frac{1}{\det R(t)}\sum\limits_{j=1}^nR_{ji}^{ad}(t)\Bigl[H_j(t)
-\sum\limits_{s=1}^kB_{js}(t)U_s(-L,t)\\ \displaystyle
-\sum\limits_{s=k+1}^nC_{js}(t)U_s(L,t)
-\sum\limits_{s=1}^n\int\limits_{-L}^LD_{js}(x,t)U_s(x,t)\,dx\Bigr].
\end{array}
$$
It is easy to see that problems~(\ref{eq:1})--(\ref{eq:3})  and~(\ref{eq:integral})
are equivalent in ${\cal G}(\Omega)$.

Given $T>0$, denote
$$\Pi^T=
\{(x,t)|-L<x<L, 0<t<T\}.
$$
In~\cite{HKm} we proved that problem~(\ref{eq:1})--(\ref{eq:3}) with
smooth initial data has a unique smooth solution in $\Pi^T$, whatsoever
$T>0$. For this purpose we used the contraction mapping principle and
obtained local smooth solution. In parallel, we obtained local a priori
estimates for the latter.
To obtain global smooth solution, we used finite iteration of the
local a priori estimates.
We also derived global apriori estimates for this solution.
To prove the existence of a generalized solution to the problem under
consideration, let us rewrite just mentioned estimates from~\cite{HKm},
with taking care of the norm of $\Lambda_i$ as well as of the norms of the
elements of $R$. Notice that the assumptions imposed on
$\Lambda_i$ and $R$ here differ from those imposed in~\cite{HKm}.
To be precise, in~\cite{HKm} $\Lambda_i$ and $R_{ij}$ for all $i,j\le n$
are assumed to be, respectively, smooth and Colombeau generalized functions
locally of bounded growth.
Referring the reader to~\cite{HKm} for details, we now write down
the final a priori estimates for a global smooth solution $U$ in a suitable
for our purposes form. Set
$$
E_U(l)=\max\{|\partial_x^l U_i(x,t)|:
(x,t)\in\overline{\Pi}^T,
1\le i\le n\},
$$
$$
E_{\Lambda,max}(l_1,l_2)=
\max\{|\partial_x^{l_1}\partial_t^{l_2} \Lambda_i(x,t)|:
(x,t)\in\overline{\Pi}^T,
1\le i\le n\},
$$
$$
E_{\Lambda,min}=
\min\{|\Lambda_i(x,t)|:
(x,t)\in\overline{\Pi}^T,
1\le i\le n\},
$$
$$
E_R=
\max\limits_{t\in[0,T]}\biggl|\frac{1}{\det R(t)}\biggr|,
$$
$$
E_B(l)=
\max\{|B_{ij}^{(l)}(t)|:
t\in[0,T],
1\le i,j\le n\},
$$
$$
E_D(l)=
\max\{|\partial_t^l D_{ij}(x,t)|:
(x,t)\in\overline{\Pi}^T,
1\le i,j\le n\},
$$
$$
E_F=\max\{|\nabla_Uf_i(x,t,y)|\,:\,(x,t,y)\in\overline{\Pi}^T\times{\bf R},1\le i\le n\}
$$
$$
\begin{array}{ccc}
\displaystyle
q_0=n^2\max\limits_{\begin{array}{cc}\scriptstyle
t\in[0,T],\\\scriptstyle
1\le i,j\le n\end{array}}\biggl
|\frac{R_{ji}^{ad}(t)}{R(t)}\biggr|
\biggl[nE_F\biggl(\max\limits_{
\begin{array}{cc}\scriptstyle
t\in[0,T],1\le j\le n,\\\scriptstyle
1\le s\le k,k+1\le r\le n
\end{array}}
\{|B_{js}(t)|,|C_{jr}(t)|\}
+2LE_{D}(0)
|\biggr)\\[7mm]
\displaystyle
+E_{D}(0)
E_{\Lambda,max}(0,0)\biggr]+nE_F,
\end{array}
$$
$$
q_m=(q_0-nE_F)E_{\Lambda,max}(0,0)^m(E_{\Lambda,min})^{-m}
+nE_F+mE_{\Lambda,max}(1,0).
$$
With this notation, we have
\begin{equation}\label{eq:25}
\begin{array}{c}
\displaystyle
E_U(m)
\le
P_{1,m}\biggl(\frac{1}{1-q_mt(m)},
\max\{E_B^n(0),E_C^n(0)\},E_R,E_D(0),
n^2,L,\\[5mm]
\displaystyle
(E_{\Lambda,max}(0,0))^m,(E_{\Lambda,min})^{-m}\biggr)
\\[5mm]
\displaystyle
\times
P_{2,m}\biggl(\max\limits_{0\le l\le m-1}E_U(l),
(E_{\Lambda,min})^{-\mbox{\rm sgn}(m)},
\max\limits_{1\le l_1+l_2\le m}E_{\Lambda,max}(l_1,l_2),
\max\limits_{0\le l\le m}\{E_B(l),E_C(l)\},
\\[5mm]
\displaystyle
\max\limits_{0\le l\le m}E_D(l),E_R,
\max\limits_{\begin{array}{cc}\scriptstyle
t\in[0,T],\\\scriptstyle
1\le i\le n,0\le l\le m\end{array}}|H_i^{(l)}(t)|,
\max\limits_{\begin{array}{cc}\scriptstyle
x\in[-L,L],\\\scriptstyle
1\le i\le n,0\le l\le m\end{array}}|A_i^{(l)}(x)|
\biggr),
\end{array}
\end{equation}
where
$t(m)\le\min\{L/E_{\Lambda,max}(0,0),
1/q_m\}$, $P_{1,m}$ is a polynomial of degree $8\lceil T/t(m)\rceil$
with positive constant coefficients not depending on $\varepsilon$, and
$P_{2,m}$ is a polynomial whose degree
depends on $m$ but neither on $T$ nor on $t(m)$
(and, therefore, not depending on $\varepsilon$)
with positive constant
coefficients depending on~$f$ and not depending on~$\varepsilon$.

The classical smooth solution to~problem~(\ref{eq:1})--(\ref{eq:3})
satisfying estimates~(\ref{eq:25}) in $\Pi^T$  for all $m\in{\bf N}_0$
can be constructed
by the sequential approximation method. This solution will serve to
build up a representative of the Colombeau solution. We now construct
the latter. Accordingly to the assumptions of the theorem, we consider all
the initial data as elements of the corresponding Colombeau algebras.
We choose representatives $\lambda$, $a$, $b$, $c$, $d$, and $h$
of $\Lambda$, $A$, $B$, $C$, $D$, and $H$, respectively, with the properties required in the
theorem. Hence a representative of $R$  is therewith defined.
We will denote it by $r$.
Let $\phi=\varphi\otimes\varphi\in {\cal A}_0({\bf R}^2)$.
Consider a prospective representative  $u=u(\phi,x,t)$ of $U$
which is the classical smooth solution to problem~(\ref{eq:1})--(\ref{eq:3})
with initial data $\lambda(\phi,x,t)$, $a(\varphi,x)$ and boundary data $b(\varphi,t)$,
$c(\varphi,t)$, $d(\phi,x,t)$, $h(\varphi,t)$. It remains  to
show the moderatness of $u$, i.e. that $u\in{\cal E}_{\mbox{\rm M}}$.
To do so, we will obtain moderate growth estimates
of $u(\phi_{\varepsilon},x,t)$ in terms of the
regularization parameter~$\varepsilon$.

Let $\varepsilon$ be small enough and
 $\phi\in{\cal A}_N({\bf R}^2)$ with $N$ chosen so large that the following conditions are
 true:

a) the moderation property holds for
$a(\varphi_{\varepsilon},x)$ and $h(\varphi_{\varepsilon},t)$;

b) the local-$\gamma$-invertibility
estimate (see Definition~\ref{defn:gainvert})
holds for $\lambda_i(\phi_{\varepsilon},x,t)$
and $r(\varphi_{\varepsilon},t)$.

c) the local-$\gamma$-growth
estimate (see Definition~\ref{defn:ga} ) holds
for $\lambda_i(\phi_{\varepsilon},x,t)$,
$b_{ij}(\varphi_{\varepsilon},t)$, $c_{ij}(\varphi_{\varepsilon},t)$
and $d_{ij}(\phi_{\varepsilon},x,t)$,
where $i\le n$ and $j\le n$.

d) the local-$\gamma_1$-growth
estimate holds
for
$\partial_x\lambda_i(\phi_{\varepsilon},x,t)$,
where $i\le n$.

It suffices to prove the moderateness of $P_{1,m}$ and $P_{2,m}$
for all $m\in{\bf N}_0$, where $U(x,t)$, $\Lambda(x,t)$,  $A(x)$, $B(t)$, $C(t)$,
$R(t)$, $D(x,t)$, and $H(t)$ are
replaced by their representatives $u(\phi,x,t)$,
$\lambda(\phi,x,t)$, $a(\varphi,x)$,  $b(\varphi,t)$,
$c(\varphi,t)$, $r(\varphi,t)$, $d(\phi,x,t)$, and $h(\varphi,t)$,
respectively.
 We see at once that for each $m\in{\bf N}_0$ the estimate
$$
q_m\le \gamma^{2N(m+1)+1}(\varepsilon)+\gamma_1^{N+1}(\varepsilon)
$$
is true for all sufficiently small  $\varepsilon$.
Since $t(m)\le\min\{L/
E_{\Lambda,max}(0,0),1/q_m\}$
and $E_{\Lambda,max}(0,0)$
$\ge 1/\gamma^{N}(\varepsilon)$ for all $\varphi\in{\cal A}_N({\bf R})$,
we can choose
$t(0)=1/[2(\gamma^{2N(m+1)+1}(\varepsilon)+\gamma_1^{N+1}(\varepsilon))]<1/q_0$.
Taking into account~(\ref{eq:ga}), for each $m\in{\bf N}_0$ and for all small enough $\varepsilon$ we have
$$
\begin{array}{ccccc}
\displaystyle
\biggl(\frac{1}{1-q_mt(m)}\biggr)^{\lceil T/t(m)\rceil}\le
2^{\lceil 2T(\gamma^{2N(m+1)+1}(\varepsilon)+\gamma_1^{N+1}(\varepsilon))\rceil}
\nonumber\\[6mm]
\displaystyle
\le
\Bigl(\gamma(\varepsilon)^{\gamma^{2N(m+1)+1}(\varepsilon)}\Bigr)
^{\lceil 2T\rceil+1}
\Bigl(\gamma_1(\varepsilon)^{\gamma_1^{N+1}(\varepsilon)}\Bigr)
^{\lceil 2T\rceil+1}=O\biggl(
\frac{1}{\varepsilon}\biggr),\nonumber\\[6mm]
\displaystyle
\biggl(
\max\{E_b^n(0),E_c^n(0)\}E_rE_d(0)n^2L
(E_{\lambda,max}(0,0))^m(E_{\lambda,min}(0))^{-m}
\biggr)^{\lceil T/t(m)\rceil}\nonumber\\[6mm]
\displaystyle
\le
\gamma(\varepsilon)^{N(2m+n+2)\lceil
2T(\gamma^{2N(m+1)+1}(\varepsilon)+
\gamma_1^{N+1}(\varepsilon))\rceil}
=O\biggl(\frac{1}{\varepsilon}\biggr)\nonumber
\end{array}
$$
It follows that for each $m\in{\bf N}_0$ there exists $N\in{\bf N}$ such that for all
$\varphi\in{\cal A}_N({\bf R})$ we have
\begin{equation}\label{eq:P1m}
\begin{array}{cc}
\displaystyle
P_{1,m}\biggl(
\frac{1}{1-q_mt(m)},
\max\{E_b^n(0),E_c^n(0)\},E_r,E_d(0)
n^2,L,\\[6mm]
\displaystyle
(E_{\lambda,max}(0,0))^m,(E_{\lambda,
min})^{-m}
\biggr)
=O\biggl(\frac{1}{\varepsilon}\biggr).
\end{array}
\end{equation}
One can easily see now that for $l=0$
\begin{equation}\label{eq:dxU}
E_u(l)
=O\biggl(\frac{1}{\varepsilon}\biggr)
\end{equation}
 for all $\varphi\in{\cal A}_N({\bf R})$ with large enough $N\in{\bf N}$.
To prove similar estimates for all derivatives of $U_i$ with respect
to $x$, we use induction on $l$. Assuming~(\ref{eq:dxU})  to hold for $l\le m$,
let us show that~(\ref{eq:dxU})  is true for $l=m+1$ as well.
Indeed, let $\varepsilon$ be small enough and $\varphi\in{\cal A}_N({\bf R})$ with $N$
chosen so large that the following conditions are true:

a) the moderateness property holds for $\partial^sa(\varphi_{\varepsilon},x)$,
$\partial^sh(\varphi_{\varepsilon},t)$,
$b^{(s)}(\varphi_{\varepsilon},t)$, $c^{(s)}(\varphi_{\varepsilon},t)$,
$\partial_t^sd(\phi_{\varepsilon},x,t)$,
$\partial_x^{l_1}\partial_t^{l_2}\lambda(\phi_{\varepsilon},x,t)$,
$\partial_x^{l}u(\phi_{\varepsilon},x,t)$
for all $0\le s\le m+1$, $0\le l\le m$,
$0\le l_1+l_2\le m+1$;

b) the local-$\gamma$-invertibility
estimate holds for $\lambda_i(\phi_{\varepsilon},x,t)$.\\
Note that $\partial_x^{l}u(\phi_{\varepsilon},x,t)$ for
$0\le l\le m$ has moderateness property due to the induction assumption.
Since $P_{2,m}$ is a polynomial whose degree does not depend on $\varepsilon$,
the moderateness of $P_{2,m}$ becomes obvious.
We are done by~(\ref{eq:P1m}).

The moderate estimates on $t$ as well as on mixed derivatives
follow immediately from~(\ref{eq:1})
by successive differentiation.
This finishes the existence part of
the proof.

The proof of the uniqueness part  follows the same scheme.
The only difference is that now we
consider problem~(\ref{eq:1})--(\ref{eq:3})
with right hand sides of~(\ref{eq:2})
and~(\ref{eq:3}) in ${\cal N}$. The analysis is even simplier since
by~\cite{gro},
 it is sufficient to check negligibility at order zero.
 The proof is complete.
\end{proof}

\begin{rem}\rm
To prove the theorem, we used an integral-operator form~(\ref{eq:integral})
of the problem under consideration. Considering~(\ref{eq:integral})
with  respect
to a Colombeau function $U\in{\cal G}(\overline\Pi^T)$, we see that the right
hand side of~(\ref{eq:integral}) includes compositions of generalized
functions. Specifically, we have compositions of the singular initial
and boundary data as well as the function $U$ with the singular
characteristic curves.

Note that the Colombeau algebra ${\cal G}$ is invariant under superposition with
smooth polynomially bounded maps. In spite of the fact that the latter
 is not the case for the compositions involved by~(\ref{eq:integral}),
all terms in~(\ref{eq:integral}) are well-defined in the Colombeau sense.
To show this,
consider system~(\ref{eq:integral})
with $U$ replaced by $u(\phi,x,t)$, where the latter is a representative of
the Colombeau solution stated in the theorem.
{}From the proof it follows that,
given $(x,t)\in\overline\Pi^T$,
the domain of dependence for $u(\phi_{\varepsilon},x,t)$
is included in a compact subset of $\overline\Pi^T$
which is independent of $\varepsilon>0$ and $\varphi\in{\cal A}_0({\bf R})$.
This means that we here do not have
the effect of infinite propagation
speed (which could be caused by
the fact that characteristic curves depend on $\varphi\in{\cal A}_0({\bf R})$).

We conclude that  the right
hand side of~(\ref{eq:integral}) is  well-defined in the Colombeau sense.
\end{rem}

\end{document}